\documentclass[12]{amsart}
\usepackage{latexsym}
\usepackage{amsthm}
\usepackage{amssymb}
\usepackage{amsfonts, amsmath}
\newtheorem{Prop}{Proposition}

\newtheorem{Lm}[Prop]{Lemma}
\newtheorem{Thm}{Theorem}

\title{Increasing unions of Stein spaces with singularities}
\author{}
\begin{document}
\maketitle
\begin{center}
\em{Youssef Alaoui}\\
\em{y.alaoui@iav.ac.ma}\\
\end{center}

\noindent
\em{Department of Mathematics,
Hassan II Institute of Agronomy}\\ 
\em{and Veterinary Sciences,
Madinat Al Irfane, BP 6202, Rabat, 10101, Morocco,}\\
%\addresses[$\dagger$]{Department of Mathematics, Institut Agronomique et V\'et\'erianire Hassan II. B.P. 6202, Rabat-Instituts, 10101. Morocco}
%\email[Y.~ Alaoui]{comp5123ster@gmail.com}
%\keywords{Stein spaces, $q$-complete spaces, $q$-convex functions, linear sets,
%$1$-convex functions with respect to a linear set.}
\date{}
\linespread{1.3}
%\maketitle
\begin{abstract}
 We show that if $X$ is a Stein space and, if $\Omega \subset X$ is exhaustable by a sequence $\Omega_1 \subset \Omega_2 \subset \ldots \subset \Omega_n \subset \ldots$ of open Stein subsets of $X$, then $\Omega$ is Stein. This generalizes a well-known result of Behnke and Stein which is obtained for $X=\mathbb{C}^n$ and solves the union problem, one of the most classical questions in Complex Analytic Geometry. When $X$ has dimension 2, we prove that the same result follows if we assume only that $\Omega \subset \subset X$ is a domain of holomorphy in a Stein normal space. It is known, however, that if $X$ is an arbitrary complex space which is exhaustable by an increasing sequence of open Stein subsets $X_1 \subset X_2 \subset \cdots \subset X_n \subset \cdots$, it does not follow in general that $X$ is holomorphically-convex or holomorphically-separate (even if $X$ has no singularities). One can even obtain 2-dimensional complex manifolds on which all holomorphic functions are constant.
\end{abstract}
{\large Key words}: Stein spaces;
$q$-complete spaces; $q$-convex functions; strictly plurisubharmonic functions.
%\begin{document}
%\maketitle \setcounter{page}{1}\noindent
\section{Introduction }
Let $X$ be a Stein space and $D\subset X$ an open subset which is
the union of an increasing sequence of Stein open subsets of $X$.\\
\hspace*{.1in}Does it follow that $D$ is necessarily Stein ?\\
\hspace*{.1in}It is known from a classical theorem due to Behnke and Stein ~\cite{ref2} that if
$D_{1}\subset D_{2}\subset\cdots\subset D_{n}\subset\cdots$ is an increasing
sequence of Stein open sets in $\mathbb{C}^{n}$, then their union
$\displaystyle\bigcup_{j\geq 1}D_{j}$ is Stein.\\
\hspace*{.1in}In 1977, Markoe ~\cite{ref4} proved the following:\\
Let $X$ be a reduced complex space which the union of an increasing
sequence $X_{1}\subset X_{2}\subset\cdots\subset X_{n}\subset\cdots$ of Stein domains.\\
\hspace*{.1in}Then X is Stein if and only if $H^{1}(X,{\mathcal{O}}_{X})=0$.\\
\hspace*{.1in}Similarly, it is known (see ~\cite{ref6}) that in an arbitrary complex space $X$ an increasing
union of Stein spaces $(X_{n})_{n\geq 0}$ is itself Stein if $H^{1}(X,{\mathcal{O}}_{X})$ is separated.\\
\hspace*{.1in}It was shown in ~\cite{ref3} that if $(D_{j})_{j\geq 1}$ is an increasing
sequence of Stein domains in a normal Stein space $X$, then $D=\displaystyle\bigcup_{j\geq 1}D_{j}$ is a domain of holomorphy.
(i.e. for each $x\in \partial{D}$ there is $f\in O(D)$ which is not holomorphically extendable through $x$).\\
\hspace*{.1in}It was proved in ~\cite{ref7} that if $X$ is a complex space and $(D_{j})_{j\geq 1}$ is an increasing sequence
of Stein open subsets of $X$, then $D=\displaystyle\bigcup D_{j}$ is $2$-complete. We recall that a complex space $X$ is said to be
$q$-complete if there exists an exhaustion function $\phi\in C^{\infty}(X,\mathbb{R})$ which is $q$-convex on the whole space $X$,
that is every point $x\in X$ has an open neighborhood $U$ isomorphic to a closed analytic set in a domain $D\subset \mathbb{C}^{n}$
such that the restriction $\phi|_{U}$ has an extension $\tilde{\phi}\in C^{\infty}(D)$ whose Levi form $L(\tilde{\phi},z)$ has at most
$q-1$ negative or zero eingenvalues at any point $z$ of $D$.\\
\hspace*{.1in}Here we solve affirmatively the above problem in the general case.
We show that if $X$ is a Stein space and, if $\Omega$ is
an increasing sequence of Stein open subsets of $X$,
then there exists an increasing sequence $(\Omega'_{\nu})_{\nu\geq 1}$ of open subsets of $\Omega$ such that
$\Omega=\displaystyle\bigcup_{\nu\geq 1}\Omega'_{\nu}$
and there are continuous strictly psh functions $\psi''_{\nu}: \Omega'_{\nu}\rightarrow ]0,+\infty[$ with the following properties\\
(a) $\psi''_{j}> 2^{\nu+2}$ on $\Omega'_{\nu+2}\backslash \Omega'_{\nu+1}$ for every $j\geq \nu+1$.\\
(b) $(\psi''_{\nu})_{\nu\geq 1}$ is stationary on every compact subset of $\Omega$.\\
This implies that the function $\psi: \Omega\rightarrow \mathbb{R}$ defined by $\psi=lim\psi''_{\nu}$
is a continuous strictly psh exhaustion function on $\Omega$.
\section{The union problem}
In order to solve the problem in dimension $2$, it is sufficient to show
\begin{Thm}
Every domain of holomorphy $D$ which is relatively compact in a $2$-dimensional normal Stein space $X$
is Stein.
\end{Thm}
\begin{proof}
By the theorem of Andreotti-Narasimhan ~\cite{ref1} we have only to prove that
$D$ is locally Stein and, we may of course assume that $X$ is connected.\\
\hspace*{.1in}Let $p\in \partial{D}\cap Sing(X)$, and choose a connected Stein open neighborhood
$U$ of $p$ with $U\cap Sing(X)=\{p\}$ and such that $U$ is biholomorphic to a closed analytic set
in a domain $M$ in some $\mathbb{C}^{N}$. Let $E$ be a complex affine subspace of $\mathbb{C}^{N}$
of maximal dimension such that $p$ is an isolated point of $E\cap U$.\\
By a coordinate transformation, one can obtain that $z_{i}(p)=0$ for all $i\in\{1,2,\cdots,N\}$
and we may assume that there is a connected Stein open neighborhood $V$ of $p$ in $M$ such that
$U\cap V\cap\{z_{1}(x)=z_{2}(x)=0\}=\{p\}$.\\
We may suppose that $N \geq 4$ and, let\\
$E_{1}=V\cap \{z_{2}(x)=\cdots=z_{N-1}(x)=0\}$, $E_{2}=\{x\in E_{1}: z_{1}(x)=0\}$.\\
Then $A=(U\cap V)\cup E_{1}$ is a Stein closed analytic set in $V$ as the union of two closed analytic subsets
of $V$.\\
Let $\xi : \tilde{A}\rightarrow A$ be a normalization of $A$. Then
$\xi : \tilde{A}\backslash \xi^{-1}(p)\rightarrow A\backslash\{p\}$ is biholomorphic and,
clearly $\xi^{-1}(A\cap E_{2})=\{x\in \tilde{A} : z_{1}(\xi(x))=\cdots=z_{N-1}(\xi(x))=0\}$
is everywhere $1$-dimensional. It follows from a theorem of Simha ~\cite{ref7} that
$\tilde{A}\backslash \xi^{-1}(A\cap E_{2})$ is Stein. Hence
$A\backslash E_{2}=\xi(\tilde{A}\backslash \xi^{-1}(A\cap E_{2}))$ itself is Stein.\\
Since $p\in E_{2}$ is the unique singular point of $A$, then $U\cap V\cap D$ is Stein, being a domain
of holomorphy in the Stein manifold $A\backslash E_{2}$.
\end{proof}
Let now $X$ be a Stein space of dimension $n\geq 2$ and $\Omega\subset X$ an open subset which is
the union of an increasing sequence $\Omega_{1}\subset \Omega_{2}\subset\cdots\subset \Omega_{n}\subset\cdots$
of Stein open sets in $X$. Let $\phi_{\nu}: \Omega_{\nu}\rightarrow ]0, +\infty[$
be a smooth strictly psh exhaustion function on $\Omega_{\nu}$,
and let $(d_{\nu})_{\nu\geq 1}$
be a sequence with $d_{\nu}<d_{\nu+1}$, and $Sup d_{\nu}=+\infty$. One may assume that
if $\Omega'_{\nu}=\{x\in \Omega_{\nu}: \phi_{\nu}(x)<d_{\nu}\}$,
then $\Omega'_{\nu}\subset\subset \Omega'_{\nu+1}$.
\begin{Lm}{-There exist for each $\nu\geq 1$ an exhaustion function $\varphi_{\nu}\in C^{\infty}(\Omega_{\nu})$
which is strictly psh in a neighborhood
of $\overline{\Omega'}_{\nu}\setminus\Omega'_{\nu-1}$, a locally finite covering $(U_{\nu})_{\nu\geq 1}$
of $\Omega$ by open sets $U_{\nu}\subset \Omega'_{\nu+1}$, and constants $c_{\nu}\in \mathbb{R}$, $\nu\geq 1$,
with the following properties:\\
(a) For each $\nu\geq 1$ there exists a function $\psi_{\nu}: \Omega'_{\nu+1}\rightarrow ]0,+\infty[$ such that $\psi_{\nu}|_{U_{\nu}}$
is strictly psh
and $\psi_{\nu}=\psi_{\nu-1}$ on $\{x\in U_{\nu}: \varphi_{\nu+1}(x)<c_{\nu}\}\cap U_{\nu-1}$.\\
(b) For every index $\nu\geq 1$, there exists $\varepsilon_{\nu}>0$ such that\\
$\Omega'_{\nu-1}\setminus\overline{\Omega'}_{\nu-2}\subset  \{x\in U_{\nu}: \varphi_{\nu+1}(x)<c_{\nu}-\varepsilon_{\nu}\}$
and\\ $\{x\in U_{\nu}: \varphi_{\nu+1}(x)<c_{\nu}+\varepsilon_{\nu}\}\subset U_{\nu-1}$.}
\end{Lm}
\begin{proof}
There exists a $C^{\infty}$ exhaustion function $\varphi_{\nu+1}$ on $\Omega_{\nu+1}$
which is strictly plurisubharmonic in a neighborhood of $\overline{\Omega'}_{\nu+1}\setminus\Omega'_{\nu}$ such that, if
$m_{\nu+1}=Min_{\overline{\Omega'}_{\nu+1}\backslash\Omega'_{\nu}} \varphi_{\nu+1}$
and $M_{\nu+1}=Max_{\overline{\Omega'}_{\nu-1}}\varphi_{\nu+1}$,
then $m_{\nu+1}>M_{\nu+1}$.\\
\hspace*{.1in}In fact, we choose $\theta_{\nu}\in C_{0}^{\infty}(\Omega_{\nu+1})$ with compact support
in $\Omega_{\nu+1}\setminus\overline{\Omega'}_{\nu-1}$ so that $0\leq \theta_{\nu}\leq 1$ and $\theta_{\nu}(x)=1$
when $x\in \overline{\Omega'}_{\nu+1}\setminus\Omega'_{\nu}$.
Let $\xi$ be a point of $\partial{\Omega'}_{\nu-1}$ such that
$\phi_{\nu+1}(\xi)=Max_{\overline{\Omega'}_{\nu-1}}\phi_{\nu+1}$.
Then it is clear that
$$\varphi_{\nu+1}=\phi_{\nu+1}+\phi_{\nu+1}(\xi)\theta_{\nu}$$
satisfies the requirements.\\
\hspace*{.1in}We now assume that $\Omega_{0}=\emptyset$ and put
$$U_{1}=\Omega'_{2}, \ \ and \ \ U_{\nu}=(\Omega'_{\nu+1}\setminus\overline{\Omega'}_{\nu-2})
 \ \ for \ \  \nu\geq 2.$$
 Then $(U_{\nu})_{\nu\geq 1}$ is a locally finite covering of $\Omega$.
Moreover, if we set
$$c'_{\nu}=m_{\nu+1}=Inf\{\varphi_{\nu+1}(x), x\in (\overline{\Omega'}_{\nu+1}\setminus\Omega'_{\nu})\},$$
then
$$(\overline{\Omega'}_{\nu-1}\setminus\overline{\Omega'}_{\nu-2})\subset \{x\in U_{\nu}: \varphi_{\nu+1}(x)<c'_{\nu}\}\subset
(\Omega'_{\nu}\setminus\overline{\Omega'}_{\nu-2})\subset U_{\nu-1}.$$
\hspace*{.1in}Furthermore, there exist $c_{\nu}>0$ and $\varepsilon_{\nu}>0$ such that\\
$c_{\nu}+\varepsilon_{\nu}=c'_{\nu}$ and
$(\overline{\Omega'}_{\nu-1}\setminus\overline{\Omega'}_{\nu-2})\subset \{x\in U_{\nu}: \varphi_{\nu+1}(x)<c_{\nu}-\varepsilon_{\nu}\}$.\\
\hspace*{.1in}Moreover, if the function $\theta_{\nu}\in C^{\infty}_{0}(\Omega_{\nu+1}\backslash\Omega'_{\nu-1})$ is chosen so that
$\theta_{\nu}=1$ on
$$(\Omega'_{\nu+1}\backslash\Omega'_{\nu})\cup\{x\in \overline{\Omega'}_{\nu}\backslash{\Omega'}_{\nu-1} :
Inf_{{\Omega'}_{\nu+1}\backslash{\Omega'_{\nu}}}\phi_{\nu+1}-\frac{\varepsilon_{\nu}}{2}\leq \phi_{\nu+1}(x)\leq Inf_{{\Omega'}_{\nu+1}\backslash{\Omega'_{\nu}}}\phi_{\nu+1}+M_{\nu+1}\},$$
then clearly we obtain
$\{x\in U_{\nu}: c_{\nu}+\frac{\varepsilon_{\nu}}{2}\leq \varphi_{\nu+1}(x)\leq c_{\nu}+\varepsilon_{\nu}\}\subset \{\theta_{\nu}=1\}$.
Therefore with such a choice of $\theta_{\nu}$ there exists for each $\nu$ a function $\psi_{\nu}: \Omega'_{\nu+1}\rightarrow ]0,+\infty[$
such that $\psi_{\nu}|_{U_{\nu}}$ is strictly plurisubharmonic and,
$\psi_{\nu}=\psi_{\nu-1}$ on $\{x\in U_{\nu}: \varphi_{\nu+1}(x)<c_{\nu}+\frac{\varepsilon_{\nu}}{2}\}.$\\
\hspace*{.1in}In fact, if $\nu=1$, then it is obvious that $\psi_{1}=\phi_{2}$ has the required properties for $\Omega_{1}=\emptyset$,
since $U_{1}=\Omega'_{2}$ and $\{x\in U_{1}: \varphi_{2}(x)<c_{1}+\frac{\varepsilon_{1}}{2}\}$ is contained in $\Omega'_{1}$.\\
\hspace*{.1in}We now assume that $\nu\geq 2$ and, that $\psi_{1}, \cdots, \psi_{\nu-1}$
have been constructed. let $\chi_{\nu}(t)=a_{\nu}(t-c_{\nu}-\frac{\varepsilon_{\nu}}{2})$
where $a_{\nu}$ is a positive constant, and consider the function
$\psi_{\nu}: \Omega'_{\nu+1}\rightarrow ]0, +\infty[$ defined by
\[\psi_{\nu}=\left\{\begin{array}{cc}
\psi_{\nu-1} \ \ on \ \  \{\varphi_{\nu+1}\leq c_{\nu}-\varepsilon_{\nu}\}\\
Max(\psi_{\nu-1}, \chi_{\nu}(\varphi_{\nu+1})) \ \ on \ \  \{c_{\nu}-\varepsilon_{\nu}\leq \varphi_{\nu+1}\leq c_{\nu}+\varepsilon_{\nu}\}\\
\chi_{\nu}(\phi_{\nu+1}+\phi_{\nu+1}(\xi)) \ \ on \ \  \{\varphi_{\nu+1}\geq c_{\nu}+\varepsilon_{\nu}\}
\end{array}
\right.\]
Since on $U'_{\nu}=\{x\in U_{\nu}: \varphi_{\nu+1}(x)<c_{\nu}+\frac{\varepsilon_{\nu}}{2}\}\subset U_{\nu-1}$ we have
$\psi_{\nu-1}>0>\chi_{\nu}(\varphi_{\nu+1})$ and
$\psi_{\nu-1}$ is strictly psh on $U_{\nu-1}$,
then $\psi_{\nu}|_{U'_{\nu}}=\psi_{\nu-1}|_{U'_{\nu}}$ is strictly psh on $U'_{\nu}$.
On the other hand, the subset $\{c_{\nu}+\frac{\varepsilon_{\nu}}{2}\leq \varphi_{\nu+1}\leq c_{\nu}+\varepsilon_{\nu}\}\subset U_{\nu-1}$
is contained in $\{\theta_{\nu}=1\}$, which implies that
$\psi_{\nu}=Max(\psi_{\nu-1}, \chi_{\nu}(\phi_{\nu+1}+\phi_{\nu+1}(\xi)))$ on
$\{c_{\nu}+\frac{\varepsilon_{\nu}}{2}\leq \varphi_{\nu+1}\leq c_{\nu}+\varepsilon_{\nu}\}$.
Then clearly the function $\psi_{\nu}$ is well-defined and satisfies the required conditions,
if $a_{\nu}$ is taken so that
$a_{\nu}\frac{\varepsilon_{\nu}}{2}>Max_{\{\varphi_{\nu+1}=c_{\nu}+\varepsilon_{\nu}\}\cap \Omega'_{\nu}}\psi_{\nu-1}$.
\end{proof}
\begin{Thm}If $X$ is a Stein space and $\Omega$
an open subset of $X$ which is an increasing union of Stein open
sets in $X$, then $\Omega$ is Stein.
\end{Thm}
\begin{proof}
We shall prove that there exists for each $\nu\geq 1$ a continuous strictly psh function $\psi''_{\nu}$ in a neighborhood of $\overline{\Omega'_{\nu}}$ such that $\psi''_{j}>2^{\nu+1}$ on $\Omega'_{\nu+2}\backslash \Omega'_{\nu+1}$ for every
$j\geq \nu+2$ and $(\psi''_{\nu})_{\nu\geq 1}$ is stationary on
every compact set in $\Omega$.\\
\hspace*{.1in}In fact, let $\varphi'_{\nu}$ be the function defined by
\[\varphi'_{\nu}=\left\{\begin{array}{cc}
\psi_{\nu} \ \ on \ \  \Omega'_{\nu+1}\backslash\overline{\Omega'}_{\nu-1}\\
\psi_{\mu} \ \ on \ \  \{x\in U_{\mu+1}: \varphi_{\mu+2}(x)< c_{\mu+1}-\varepsilon_{\mu+1}\} \ \ for \ \ \mu\leq \nu
\end{array}
\right.\]
\hspace*{.1in}Then, by lemma $1$, $\varphi'_{\nu}$ is a continuous strictly plurisubharmonic function on $\Omega'_{\nu+1}$.\\
\hspace*{.1in}Moreover, we have $\varphi'_{\nu}=\varphi'_{\nu-1}$ on $\{x\in U_{\mu+1}: \varphi_{\mu+2}(x)<c_{\mu+1}-\varepsilon_{\mu+1}\}$ for all $\mu\leq \nu-1$.\\
\hspace*{.1in}Let now $K$ be a compact set in $\Omega$ and $\nu\geq 2$ such that
$K\subset \Omega'_{\nu-1}$. Since $\varphi'_{\nu}=\varphi'_{\nu-1}$ on
$K\cap (\overline{\Omega'}_{\mu}\setminus\overline{\Omega'}_{\mu-1})\subset \{x\in U_{\mu+1}: \varphi_{\mu+2}(x)<c_{\mu+1}-\varepsilon_{\mu+1}\}$
for all $\mu\leq \nu-1$, then $\varphi'_{\nu}=\varphi'_{\nu-1}$ on $K.$
This implies that the sequence $(\varphi'_{\nu})_{\nu\geq 1}$ is stationary on every compact subset of $\Omega$.\\
\hspace*{.1in}Let now $\nu\geq 1$ be an arbitrary natural number. Then there exists
a smooth function $\psi'_{\nu}\in C^{\infty}(X)$ which is strictly plurisubharmonic in
a neighborhood of $(X\backslash\Omega'_{\nu+1})\cup \overline{\Omega'}_{\nu}$ such that
$\psi'_{\nu}> 2^{\nu+2}$ in $\overline{\Omega'}_{\nu+2}\backslash \Omega'_{\nu+1}$ but
$\psi'_{\nu}<0$ in $\overline{\Omega'}_{\nu}$.\\
\hspace*{.1in} In fact, let $h\in C^{\infty}(X)$ be a strictly plurisubharmonic exhaustion function
such that $h<0$ in $\overline{\Omega'}_{\nu}$, and
let $\chi_{\nu}\in C^{\infty}(X)$ be a smooth function with compact support
in $\Omega'_{\nu+1}$ such that $\chi_{\nu}=1$ in $\overline{\Omega'}_{\nu}$. Then it is clear that
$$h_{\nu}=h+b_{\nu}.\chi_{\nu},$$
where $b_{\nu}=Min_{x\in \overline{\Omega'}_{\nu+2}\backslash \Omega'_{\nu+1}} h(x)$,
is a smooth exhaustion function on $X$ which is strictly plurisubharmonic in a neighborhood
of $(X\backslash\Omega'_{\nu+1})\cup \overline{\Omega'}_{\nu}$ such that if
$m'_{\nu}=Min_{y\in \overline{\Omega'}_{\nu+2}\backslash\Omega'_{\nu+1}} h_{\nu}(y)$
and $M'_{\nu}=Max_{y\in \overline{\Omega'}_{\nu}}h_{\nu}(y)$, then $m'_{\nu}>M'_{\nu}$.\\
Let $\varepsilon'_{\nu}>0$ be such that $m'_{\nu}>M'_{\nu}+\varepsilon'_{\nu}$. Then we can choose a sufficiently big constant $C_{\nu}>1$
so that
$$\psi'_{\nu}(x)=C_{\nu}(h_{\nu}(x)-M'_{\nu}-\varepsilon'_{\nu})$$
is $>2^{\nu+2}$ in $\overline{\Omega'}_{\nu+2}\backslash\Omega'_{\nu+1}$, $\psi'_{\nu}<0$ in $\overline{\Omega'}_{\nu}$,
and strictly plurisubharmonic in a neighborhood of $(X\backslash\Omega'_{\nu+1})\cup\overline{\Omega'}_{\nu}$.\\
\hspace*{.1in}If now we consider the following function defined in lemma $1$
\[\psi_{\nu}=\left\{\begin{array}{cc}
\psi_{\nu-1} \ \ on \ \  \{\varphi_{\nu+1}\leq c_{\nu}-\varepsilon_{\nu}\}\\
Max(\psi_{\nu-1}, \chi_{\nu}(\varphi_{\nu+1})) \ \ on \ \  \{c_{\nu}-\varepsilon_{\nu}\leq \varphi_{\nu+1}\leq c_{\nu}+\varepsilon_{\nu}\}\\
\chi_{\nu}(\phi_{\nu+1}+\phi_{\nu+1}(\xi)) \ \ on \ \  \{\varphi_{\nu+1}\geq c_{\nu}+\varepsilon_{\nu}\}
\end{array}
\right.\]
and the fact that $c_{\nu}+\varepsilon_{\nu}=Inf\{\varphi_{\nu+1}(x), x\in \overline{\Omega'}_{\nu+1}\backslash\Omega'_{\nu}\}$, we find that
$$\Omega'_{\nu+1}\backslash\overline{\Omega'}_{\nu}\subset\{x\in U_{\nu} : \varphi_{\nu+1}(x)\geq c_{\nu}+\varepsilon_{\nu}\}$$
and, on the set $\Omega'_{\nu+1}\backslash\overline{\Omega'}_{\nu}$ we have
$$\varphi'_{\nu}=\psi_{\nu}=\chi_{\nu}(\phi_{\nu+1}+\phi_{\nu+1}(\xi))\geq a_{\nu}(\varphi_{\nu+1}-c_{\nu}-\frac{\varepsilon_{\nu}}{2})\geq a_{\nu}\frac{\varepsilon_{\nu}}{2}.$$
\noindent
We can therefore choose $a_{\nu}$ again big enough so that
$a_{\nu}.\frac{\varepsilon_{\nu}}{2}>\psi'_{\nu}$ on $\overline{\Omega'}_{\nu+1}\backslash\Omega'_{\nu}$.
Moreover, by suitable choice of the constants $a_{\mu}$ we can also achieve that $\varphi'_{\nu}>\psi'_{\mu}$ on $\Omega'_{\mu+1}\backslash\Omega'_{\mu}$ for
all $\mu<\nu$.
In fact, since $\Omega'_{\mu}\backslash\overline{\Omega'}_{\mu-1}\subset \{x\in U_{\mu+1} : \varphi_{\mu+2}(x)<c_{\mu+1}-\varepsilon_{\mu+1}\}$,
then, for every $2\leq \mu\leq \nu$, $\varphi'_{\nu}=\psi_{\mu}$ on $\Omega'_{\mu}\backslash\overline{\Omega'}_{\mu-1}$. If
we set
$A_{\mu}=\Omega'_{\mu}\backslash\overline{\Omega'}_{\mu-1}\cap \{x\in U_{\mu} : \varphi_{\mu+1}(x)<c_{\mu}-\varepsilon_{\mu}\},$
then $\psi_{\mu}=\psi_{\mu-1}$ on $A_{\mu}$. Since in addition $\Omega'_{\mu}\backslash\overline{\Omega'}_{\mu-1}\subset \{x\in U_{\mu-1} : \varphi_{\mu}(x)\geq c_{\mu-1}+\varepsilon_{\mu-1}\}$, then on the set $A_{\mu}$ we have
$\varphi'_{\nu}=\psi_{\mu}=\psi_{\mu-1}\geq \chi_{\mu-1}(\varphi_{\mu})\geq a_{\mu-1}\frac{\varepsilon_{\mu-1}}{2}$.
Let now $x\in \Omega'_{\mu}\backslash \Omega'_{\mu-1}$. If $x\notin A_{\mu}$, since $x\in U_{\mu}$, then $\varphi_{\mu+1}(x)\geq c_{\mu}-\varepsilon_{\mu}$.
Because
$\Omega'_{\mu}\backslash\overline{\Omega'}_{\mu-1}\subset \{x\in U_{\mu-1} : \varphi_{\mu}(x)\geq c_{\mu-1}+\varepsilon_{\mu-1}\},$
we obtain, if $\varphi_{\mu+1}(x)\leq c_{\mu}+\varepsilon_{\mu}$,
$$\varphi'_{\nu}(x)=\psi_{\mu}(x)=Max(\psi_{\mu-1}(x),\chi_{\mu}(\varphi_{\mu+1}(x))\geq \psi_{\mu-1}(x)=\chi_{\mu-1}(\varphi_{\mu}(x))>a_{\mu-1}\frac{\varepsilon_{\mu-1}}{2}.$$
$$\text{Or} \ \ \varphi'_{\nu}(x)=\psi_{\mu}(x)\geq \chi_{\mu}(\varphi_{\mu+1})(x)>a_{\mu}\frac{\varepsilon_{\mu}}{2}, \ \ \text{if} \ \ \varphi_{\mu+1}(x)\geq c_{\mu}+\varepsilon_{\mu}.$$
So we may of course take the constants $a_{\mu}$ sufficiently large so that $a_{\mu-1}.\frac{\varepsilon_{\mu-1}}{2}>\psi'_{\mu-1}$
and $a_{\mu}.\frac{\varepsilon_{\mu}}{2}>\psi'_{\mu-1}$
on $\overline{\Omega'}_{\mu}\backslash\Omega'_{\mu-1}$ for all $\mu\leq \nu$. Since only finitely many conditions are required
to get $\varphi'_{\nu}>\psi'_{\mu}$ on $\Omega'_{\mu+1}\backslash \Omega'_{\mu}$ for $\mu\leq \nu$,
it follows that the function $\psi''_{\nu}: \Omega'_{\nu+1}\rightarrow \mathbb{R}$ given by
$\psi''_{\nu}=Max(\varphi'_{\nu},\psi'_{\nu}, \psi'_{\nu-1}, \cdots, \psi'_{1})$
is obviously continuous and strictly plurisubharmonic in $\Omega'_{\nu+1}$.
Also it is clear that for every $j\geq \nu+1$, $\psi''_{j}\geq \psi'_{\nu}>2^{\nu+2}$ on $\Omega'_{\nu+2}\setminus \overline{\Omega'}_{\nu+1}$. \\
\hspace*{.1in}Let now
$K\subset \Omega$ be a compact subset and $\nu\geq 2$
such that $K\subset \Omega'_{\nu-1}$.
Since $\varphi'_{\nu}>0>\psi'_{\nu}$ on $\overline{\Omega'}_{\nu-1}$
and $\varphi'_{\nu}=\varphi'_{\nu-1}$ on $K$,
then
$Max(\varphi'_{\nu-1},\psi'_{\nu-1},\psi'_{\nu-2},\cdots,\psi'_{1})=Max(\varphi'_{\nu},\psi'_{\nu}, \psi'_{\nu-1},\cdots, \psi'_{1})$
on $K$, which implies that the sequence $(\psi''_{\nu})_{\nu\geq 1}$
is stationary on every compact subset of $\Omega$.\\
\hspace*{.1in}This proves that the limit $\psi''$ of $(\psi''_{\nu})$ is a continuous strictly plurisubharmonic exhaustion
function on $\Omega$, which shows that $\Omega$ is Stein.
\end{proof}

\end{document}